\documentclass[a4paper,abstracton]{scrartcl}

\usepackage{type1cm}        
\usepackage{makeidx}         
\usepackage{graphicx}        
\usepackage{multicol}        
\usepackage[bottom]{footmisc}

\usepackage{newtxtext}       %
\usepackage{newtxmath}       

\usepackage{todonotes}

\makeindex             

\usepackage[utf8]{inputenc}
\usepackage{url} 
\usepackage[hidelinks]{hyperref} 
\newcommand{\email}[1]{\texttt{#1}}

\newcounter{mymac@matlab}
\newcommand{\matlab}{MATLAB%
  \ifnum\value{mymac@matlab}<1%
  \textsuperscript{\textregistered}%
  \setcounter{mymac@matlab}{1}%
  \fi%
}
\DeclareMathOperator{\dif}{d \!}
\newcommand{\mess}{\textsf{M.E.S.S.}}
\newcommand{\mmess}{\mbox{\textsf{M-}\mess{}}}
\newcommand{\morlab}{\mbox{MORLAB}}

\usepackage{authblk}
\usepackage{tikz} 
\usepackage{pgfplots}
\pgfplotsset{
/pgfplots/colormap={hot2}{[1cm]rgb255(0cm)=(0,0,0) rgb255(3cm)=(255,0,0)
rgb255(6cm)=(255,255,0) rgb255(8cm)=(255,255,255)}
}
\pgfplotsset{
colormap={redyellow}{rgb255(0cm)=(255,0,0);
rgb255(1cm)=(255,255,0)}
}
\pgfplotsset{compat = 1.13, 
  colormap name = hot2, 
  unbounded coords = jump, 
  every axis plot/.append style = {%
    line width = 1.5pt 
  } 
} 
\usepgfplotslibrary{external} 
\tikzset{external/system call = {%
    lualatex \tikzexternalcheckshellescape
    -halt-on-error
    -interaction=batchmode
    -jobname "\image" "\texsource"}} 
\tikzexternaldisable{} 

\usepackage{filemod}

\title{A Non-stationary Thermal-Block Benchmark Model for Parametric Model Order
  Reduction}
\author[1]{Stephan Rave}
\affil[1]{University of Münster, Orleans-Ring 10,
  48149 Münster, Germany, \email{stephan.rave@uni-muenster.de}}
\author[2]{Jens Saak}\affil[2]{Sandtorstr. 1, 39106 Magdeburg, Germany,
  Max Planck Institute for Dynamics of Complex Technical Systems,
  \email{saak@mpi-magdeburg.mpg.de}}
\date{February 28, 2020}
\begin{document}
\maketitle
\begin{abstract}
In this contribution we aim to satisfy the demand for a publicly
available benchmark for parametric model order reduction that is scalable
both in degrees of freedom as well as parameter dimension.
\end{abstract}

\section{Introduction}\label{sec:RS20:introduction}
Model order reduction (MOR) of parametric problems (PMOR) is accepted to be an
important field of research, in particular due to its relevance for multi-query
applications such as uncertainty quantification, inverse problems
or parameter studies in the engineering sciences. Still,
publicly available software is often either tailored to a very specific problem
or bound to a specific PDE discretization software. The joint feature of the
software packages, \texttt{emgr}~\cite{morHim18a}, \mmess{}~\cite{SaaKB-mmess-all-versions},
\morlab{}~\cite{morBenW19b} and pyMOR~\cite{morpymorweb}, reported in this volume is the attempt to make (P)MOR
available in a more general purpose fashion.
Further packages that fall into this category are
rbMIT~\cite{rbmit},
RBmatlab~\cite{rbMatlabweb,morHaa17},
RBniCS~\cite{rbnicsweb,HesthavenRozzaEtAl2016},
redbKIT~\cite{redbKIT, QuarteroniManzoniEtAl2016},
psssMOR~\cite{psssMORweb}.
So far, comparison of PMOR methods is a
difficult task~\cite{morBauBHetal17}. We think that one of the difficulties is
the lack of models that can be easily used and fairly compared in all
packages. It is the goal of this benchmark to overcome some of the shortcomings
of available benchmarks.

The MOR community Wiki~\cite{morWiki} already provides a number of parametric
benchmark models. However, most of them have either rather large dimensions
making them difficult to access directly for dense matrix-based packages
like~\cite{morBenW19b}, but also cumbersome to use during development and
testing of new sparse methods. Other benchmarks have rather limited parameter
dimension, i.e.\ they feature only scalar or at most two-dimensional
parameters. A very common feature among the benchmarks in the Wiki is that
essentially all of them are matrix-based, giving easy access for \matlab{}-based
solvers, but at the same time making it difficult for packages like
pyMOR~\cite{morpymorweb,morMilRS16,morBalMRetal19} to show their full flexibility.

Therefore, the new benchmark introduced in this chapter has a few features
addressing exactly these problems. The model is of limited dimension in the
basic version provided as matrices. On the other hand, it also provides the
FEniCS~\cite{fenics15,LogMW12} based procedural setup\footnote{Actually, the core
  feature is the unified form language (UFL)~\cite{ufl} that also other
  packages, like e.g.\ firedrake~\cite{RatHMetal16} use.} allowing for easy
generation of larger versions or integration into FEniCS-based software
packages. The current version features one to four parameters, but the setup can
be extended to higher parameter dimensions by tweaking the basic domain description
given as plain text input for gmsh~\cite{GeuR10}. Thus, we provide maximal
flexibility with a small, but scalable benchmark with  up to four
independent parameters given in a description that can easily be adapted for many PDE
discretization tools.
The benchmark we introduce here is a specific version of the so called thermal-block
benchmark. This type of model has been a standard test case in the reduced basis
community for many years, e.g.~\cite{morRozHP08}. This specific model setup is also known as the
``cookie baking problem''~\cite{morBalK16} in the numerical linear algebra
community. It further presents a flattened 2d version of what is sometimes
referred to as the ``skyscraper model'' in high performance computing,
e.g.~\cite[p.~216]{DolJN15}. We choose the common name used for this type of
model in the reduced-basis community. 

The remainder of this chapter is organized as follows. The next section provides
a basic, abstract description of the model problem. After that, in
Section~\ref{sec:RS20:variants}, we present three variants of our model that
will be used in the numerical experiments of the following chapters.
\section{Problem Description}\label{sec:RS20:description}
\begin{figure}[ht]
  \centering
  \definecolor{Gamma_in}{HTML}{ca0020}
\definecolor{Gamma_N}{HTML}{92c5de}
\definecolor{Gamma_D}{HTML}{0571b0}

\begin{tikzpicture}[x=.5\linewidth, y=.5\linewidth,thick]
  \foreach \x in {0.2,0.4,0.6,0.8} 
  {
    \draw[loosely dotted] (\x,0)--(\x,1);
    \draw[loosely dotted] (0,\x)--(1,\x);
    \draw (\x,0)--(\x,-.025);
    \draw (0,\x)--(-.025,\x);
    \node[below] at (\x, -.05) {\x};
    \node[left] at (-.05,\x) {\x};
  }
  \node[below left=.05] at (0,0) {(0,0)};
  \node[above left=.05] at (0,1) {(0,1)};
  \node[below right=.05] at (1,0) {(1,0)};
  \node[above right=.05] at (1,1) {(1,1)};
  \node at (.5,.5){\large\(\Omega_{0}\)};
  \foreach \x/\y/\n in { .3/.3/1, .7/.3/2, .7/.7/3, .3/.7/4 }
     {
     \draw[black] (\x,\y) circle [radius=.1];
     \node at (\x,\y) {\large\(\Omega_{\n}\)};
     }
  \foreach \x/\y/\n in { -.1/.5/in, .5/-.075/N, .5/1.05/N, 1.1/.5/D}
    \node[color=Gamma_\n] at (\x,\y) {\large\(\Gamma_{\n}\)};
  
  \draw[Gamma_in] (0,0)--(0,1);
  \draw[Gamma_N] (0,0)--(1,0);
  \draw[Gamma_N] (0,1)--(1,1);
  \draw[Gamma_D] (1,0)--(1,1);
  \foreach \x/\y/\n in { 0/0, 1/0, 0/1, 1/1}
    \filldraw[black] (\x,\y) circle [radius=1pt];

\end{tikzpicture}
  \caption{computational domain and boundaries}%
  \label{fig:RS20:domain}
\end{figure}
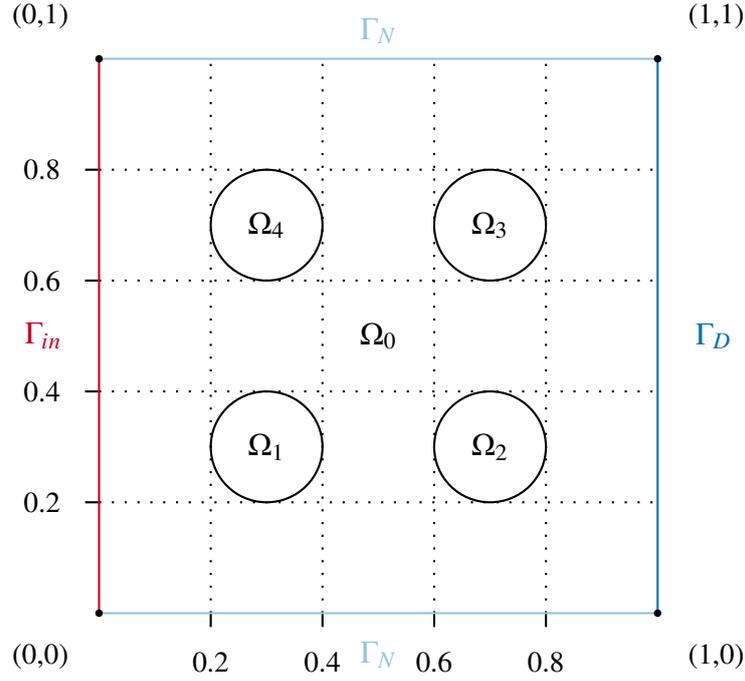
We consider a basic parabolic `thermal-block'-type benchmark problem. To this
end, consider the computational domain $\Omega := {(0, 1)}^2$ which we partition
into subdomains
\begin{align*}
    \Omega_1&:=\{\xi \in \Omega\ |\ |\xi - (0.3, 0.3)| < 0.1\},&
    \Omega_2&:=\{\xi \in \Omega\ |\ |\xi - (0.7, 0.3)| < 0.1\},\\
    \Omega_3&:=\{\xi \in \Omega\ |\ |\xi - (0.7, 0.7)| < 0.1\},&
    \Omega_4&:=\{\xi \in \Omega\ |\ |\xi - (0.3, 0.7)| < 0.1\},\\
    \Omega_0&:=\Omega \setminus (\Omega_1 \cup \Omega_2 \cup \Omega_3 \cup \Omega_4),
\end{align*}
with its boundary partitioned into
\begin{equation*}
    \Gamma_{in}:=\{0\} \times (0,1), \quad
    \Gamma_{D}:=\{1\} \times (0,1), \quad
    \Gamma_{N}:=(0,1) \times \{0,1\},
\end{equation*}
cf.~Fig.~\ref{fig:RS20:domain}.
Given a parameter $\mu \in \mathbb{R}_{\geq 0}^4$, let the heat conductivity
$\sigma(\xi; \mu)$ given by
\begin{equation}\label{eq:RS20-def_sigma}
    \sigma(\xi; \mu) :=
    \begin{cases}
        1 & \xi \in \Omega_0 \\
        \mu_i & \xi \in \Omega_i, \quad 1\leq i \leq 4,
    \end{cases}
\end{equation}
and let the temperature $\theta(t, \xi; \mu)$ in the time interval $[0,T]$ for thermal
input $u(t)$ at $\Gamma_{in}$ be given by
\begin{align*}
    \partial_t \theta(t, \xi; \mu) + \nabla \cdot (- \sigma(\xi; \mu) \nabla \theta(t, \xi; \mu)) &= 0 & t\in (0,T),&\ \xi \in \Omega,\\
    \sigma(\xi; \mu) \nabla \theta(t, \xi; \mu) \cdot n(\xi) &= u(t) & t \in (0,T),&\ \xi \in \Gamma_{in},\\
    \sigma(\xi; \mu) \nabla \theta(t, \xi; \mu) \cdot n(\xi) &= 0 & t \in (0,T),&\ \xi \in \Gamma_N,\\
    \theta(t, \xi; \mu) &= 0 & t \in (0,T),&\ \xi \in \Gamma_{D},\\
    \theta(0, \xi; \mu) &= 0 & &\ \xi \in \Omega.
\end{align*}
More precisely, we let $\theta \in L^2(0, T; V)$ with 
$\partial_t \theta(\mu) \in L^2(0, T; V^\prime)$ be given as the solution of the weak
parabolic problem
\begin{align}
    \langle \partial_t\theta(t,\cdot; \mu),\, v \rangle +
    \int_{\Omega} \sigma(\mu) \nabla \theta(t,\xi; \mu) \cdot \nabla v \dif \xi
    &= \int_{\Gamma_{in}} u(t)v \dif s & t \in (0,T),&\ v\in V,\label{eq:RS20-weak_form} \\
    \theta(0,\xi; \mu) &= 0,
\end{align}
where $V:=\{v \in H^1_0(\Omega)\ |\ v_{\Gamma_D} = 0\}$ denotes the space of Sobolev functions with vanishing trace on
$\Gamma_D$ and $V^\prime$ is its continuous dual.

As outputs $y(t; \mu) \in \mathbb{R}^4$ we consider the average temperatures in
the subdomains $\Omega_i$, i.e.
\begin{equation}\label{eq:RS20-output}
    y_i(t; \mu) := \frac{1}{|\Omega_i|} \int_{\Omega_i} \theta(t,\xi; \mu) \dif \xi, \qquad 1\leq i \leq 4.
\end{equation}

To ease the notation, we drop the explicit dependence on \(\xi\) in the following.

In view of the definition~\eqref{eq:RS20-def_sigma} of $\sigma$ as a linear combination
of characteristic functions, we can write~\eqref{eq:RS20-weak_form}--\eqref{eq:RS20-output}
as
\begin{equation}\label{eq:RS20-weak_form_abstract}
    \begin{aligned}
        \partial_t m(\theta(t; \mu), v) + a_0(\theta(t; \mu), v)
        + \sum_{i=1}^4 \mu_i \cdot a_i(\theta(t; \mu), v) &= \varphi(v)\cdot u(t) \in V \\
        \theta(0; \mu) &= 0 \\
        y_i(t; \mu) &=  \psi_i(\theta(t; \mu)), 
    \end{aligned}
\end{equation}
for $t \in (0, T)$, $v \in V$, $1 \leq i \leq 4$, with bilinear forms
$m, a_i \in \operatorname{Bil}(V, V)$
given by
\begin{equation*}
    m(w, v) := \int_\Omega wv \dif \xi \quad\text{and}\quad
    a_i(w, v) := \int_{\Omega_i} \nabla w \cdot \nabla v \dif \xi
\end{equation*}
and linear forms $\varphi, \psi_i \in V^\prime$ given by
\begin{equation*}
    \varphi(v) := \int_{\Gamma_{in}} v \dif s, \quad\text{and}\quad
    \psi_i(v) := \frac{1}{|\Omega_i|} \int_{\Omega_i} v \dif \xi.
\end{equation*}

To arrive at a discrete approximation of~\eqref{eq:RS20-weak_form_abstract},
we perform a Galerkin projection onto a space $S^1(\mathcal{T}) \cap V$
of linear finite elements w.r.t.\ a simplicial triangulation of $\Omega$
approximating the decomposition into the subdomains $\Omega_0, \ldots, \Omega_4$.
Assembling matrices $E\in \mathbb{R}^{n\times n}$ for \(m\), $A_i \in
\mathbb{R}^{n\times n}$ for \(a_{i}\),
$B \in \mathbb{R}^{n \times 1}$ for \(\varphi\), and $C \in \mathbb{R}^{4 \times n}$ for
$\psi_i$, all w.r.t.\ the finite-element basis we arrive at the linear time-invariant system
\begin{equation}
  \begin{aligned}
    E \cdot \partial_t x(t; \mu) &= A_0 \cdot x(t; \mu) + \sum_{i=1}^4 \mu_i A_i \cdot x(t; \mu) + B \cdot u(t)\\
    y(t; \mu) &= C \cdot x(t; \mu).
  \end{aligned}
  \label{eq:RS20-statespace}
\end{equation}
Here, $n$ denotes the dimension of the finite-element space and $x$ is the
coefficient vector of the discrete solution state $\theta$ w.r.t.\ the finite-element
basis.

For the numerical experiments in the following chapters, the mesh $\mathcal{T}$
was generated with gmsh version 3.0.6 with `clscale' set to $0.1$, for which the
system matrices were assembled using FEniCS 2019.1.
\begin{center}
  \fbox{%
    \parbox{.9\linewidth}{
      The source code of the model implementation as well as the resulting
      system matrices are available at   
      \begin{center} 
        \url{https://doi.org/10.5281/zenodo.3691894}
      \end{center}
    }%
  }%
\end{center}
Note that due to the handling of Dirichlet constraints in FEniCS, all matrices
were assembled over the full unconstrained space $S^1(\mathcal{T})$.
Rows of $E, A_i$ corresponding to degrees of freedom located on $\Gamma_D$
have zero off-diagonal entries.
The corresponding diagnal entries are $1$ for $E$, $-1$ for $A_0$ and
$0$ for $A_1, \ldots, A_4$.
Rows of $B$ corresponding to Dirichlet degrees of freedom are set to 0.
Consequently, all system matrices $A(\mu) := A_0 + \sum_{i=1}^4 \mu_i A_i$
have a $k$ dimensional eigenspace with eigenvalue $-1$ spanned by the
$k$ finite-element basis functions associated with $\Gamma_{D}$.

\section{Problem Variants}\label{sec:RS20:variants}
The following chapters test the model introduced in the previous section in
three different variants. The simplest case is a basic 
non-parametric version with all parameters fixed. For the parametric versions either all four
parameters are considered independent, or they are all scaled versions of a
single scalar parameter. This section introduces all of them with the specific
parameter selections and allowed parameter domains.

\begin{figure}[tbp]
  \centering
  \includegraphics[width=.8\linewidth]{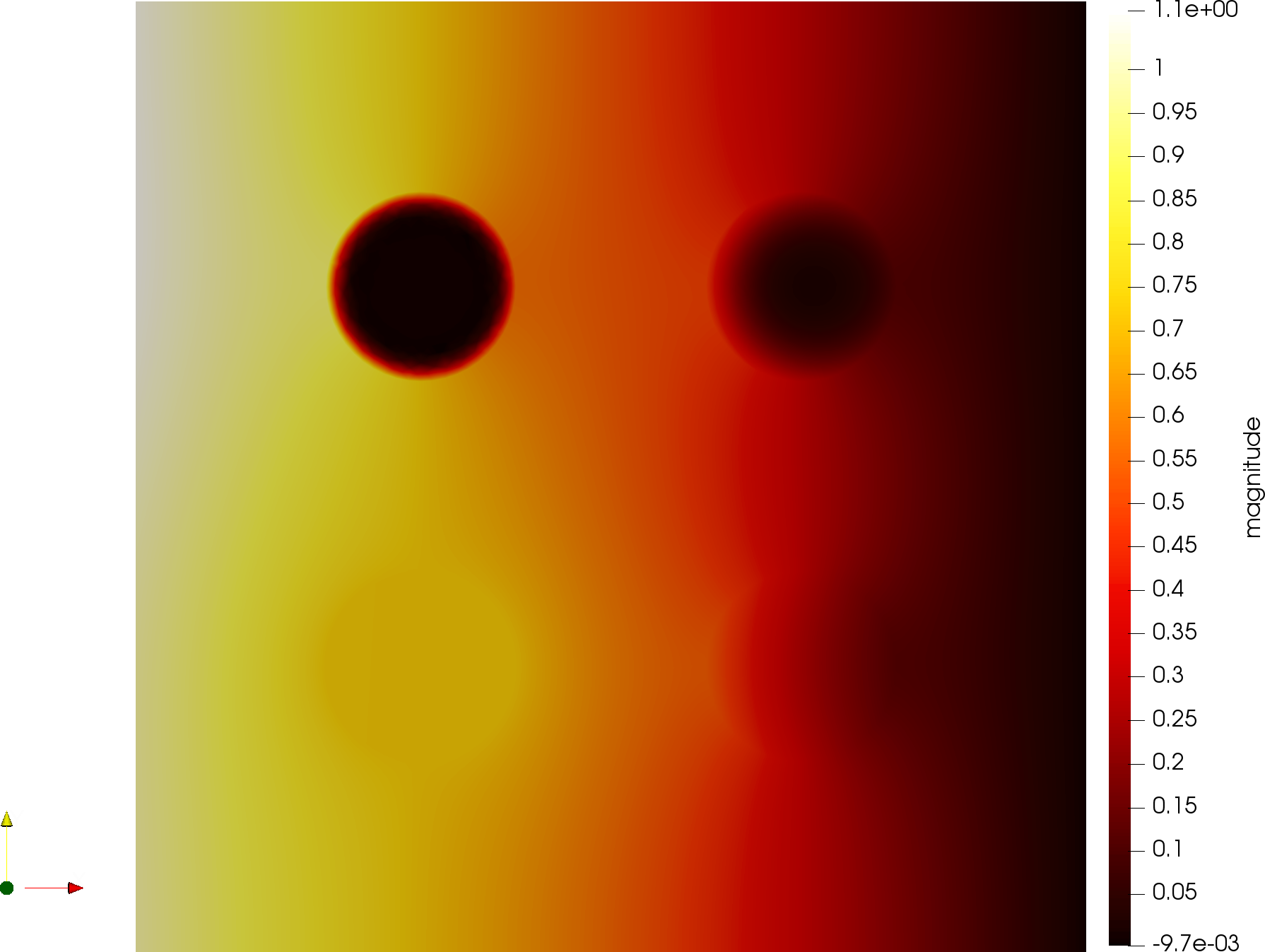}
  \caption{A sample final heat distribution.}\label{fig:RS20-final-fourp}
\end{figure}
\subsection{Four Parameter LTI System}\label{sec:RS20:quad-variant}
This represents exactly the model in~\eqref{eq:RS20-weak_form_abstract},
or~\eqref{eq:RS20-statespace}, with its full flexibility with respect to the
parameters. Note that by construction the model becomes singular in case any of
the \(\mu_{i}\) become zero.  Thus, we limit the \(\mu_{i}\) from below by
\(10^{-6}\). This will also limit the condition numbers of the linear systems
involving the matrices \(E\) and \(A_{i}\) (\(i=0,\dots,4\)) in the PDE solvers
as well as MOR routines. At the same time, we do not allow for the subdomain
heat conductivities to be drastically larger than the conductivity for
\(\Omega_{0}\). So we limit also from the above, resulting in parameter domains
\(\mu_{i}\in [10^{-6},10^{2}]\), 
(\(i=1,\dots,4\)).

Figure~\ref{fig:RS20-final-fourp} shows the final heat distribution, at \(t=1\)
after 100 steps of implicit Euler with \(\mu =[10^{2}, 10^{-2}, 10^{-3},
10^{-4}] \), in pyMOR 2019.2. 
\subsection{Single Parameter LTI System}\label{sec:RS20:single-variant}
\begin{figure}[tbp]
  \centering
  \tikzexternalenable%
  \tikzsetnextfilename{FOM_sigma}%
  \filemodCmp{graphics/FOM_sigma.tikz}{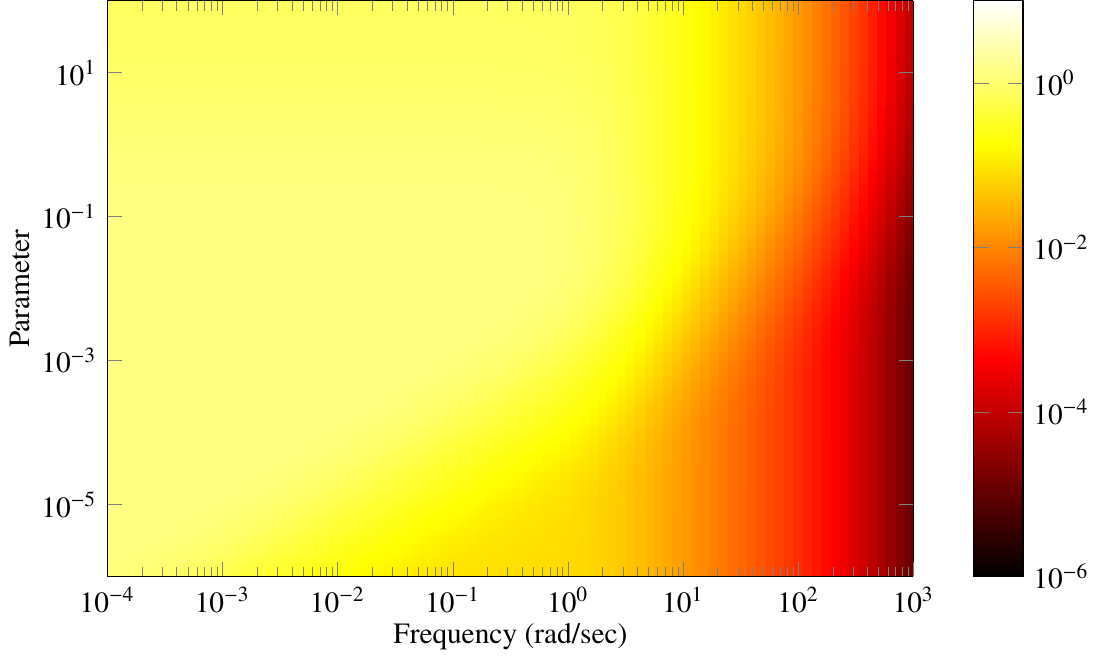}%
  {\tikzset{external/remake next}}{}%
  \begin{tikzpicture} 
  \pgfplotstableread{./figures/FOM_fom.dat}\tableROM 
  \begin{loglogaxis}[ 
    view={0}{90},
	colorbar,
	colorbar style = { 
      yticklabel = $10^{\pgfmathparse{\tick} 
        \pgfmathprintnumber\pgfmathresult}$
    },
    width = .7\textwidth, 
    height = .5\textwidth, 
    scale only axis, 
    xmin = 1e-4, 
    xmax = 1e+3, 
    ymin = 1e-6, 
    ymax = 1e+2, 
    xtick = {1e-4, 1e-3, 1e-2, 1e-1, 1e0, 1e+1, 1e+2, 1e+3}, 
    ytick = {1e-5, 1e-3, 1e-1, 1e+1}, 
    zmode = log, 
    log base z = 10, 
    point meta min = -6, 
    point meta max = 1, 
    mesh/ordering = y varies, 
    mesh/rows = 100, 
    mesh/cols = 100, 
    xlabel = {\small Frequency (rad/sec)}, 
    xlabel style = {yshift = .3em}, 
    ylabel = {\small Parameter}, 
    ylabel style = {yshift = -.3em}, 
    scaled x ticks = false, 
    x tick label style = {/pgf/number format/fixed}] 

    \addplot3[surf, shader = flat] 
    table[x index = 0, y index = 1, z index = 2] {\tableROM}; 

  \end{loglogaxis} 
\end{tikzpicture} %
  \tikzexternaldisable%

  \caption{Sigma magnitude plot for the single parameter LTI system.}%
  \label{fig:RS20-FOM-sigma}
\end{figure}
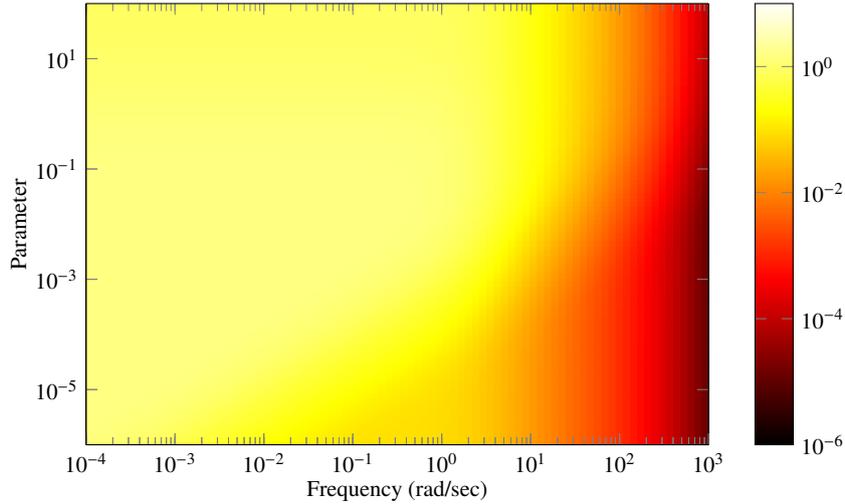
In this variation of the model the parameters are limited in flexibility. We
make them all use the same order of magnitude by defining
\begin{equation*}
  \mu = \tilde\mu\cdot
  \begin{bmatrix}
    0.2\\0.4\\0.6\\0.8
  \end{bmatrix}, 
\end{equation*}
for a single scalar parameter \(\tilde\mu\in[10^{-6},10^{2}]\). The transfer
function, arising after Laplace transformation of~\eqref{eq:RS20-statespace}
is a rational matrix-valued function of the frequency and the
parameters. Its Sigma-magnitude plot, i.e.\ the maximum singular value of the
transfer function matrix, with this restriction on \(\mu\), is shown in
Figure~\ref{fig:RS20-FOM-sigma}. 
\subsection{Non-parametric LTI System}\label{sec:RS20:lti-variant}
This is the simplest version of the benchmark. We use the setup described in
Section~\ref{sec:RS20:single-variant} with \(\tilde\mu=\sqrt{10}\). Note that
this value of \(\tilde\mu\) is rather arbitrary. Depending on the desired
application, different values may be insightful. For both time domain and
frequency domain investigations variation is strongest in the parameter range
\([10^{-5},10^{-1}]\). On the other hand values between \(1.25\) and \(5.0\)
essentially turn the model into a simple 
heat equation on the unit square with almost homogeneous heat conductivity
\(\sigma(t,\xi,\tilde\mu) \approx 1\). Hence,
\(\tilde\mu=\sqrt{10}\approx 3.1623\) appears to be a proper choice to get
reasonably close to an easy to solve textbook problem, here. Smaller values of
\(\mu\), especially when approaching \(\mu=0\), can be used to make the problem
arbitrarily ill conditioned.

\section*{Conclusion}\label{sec:RS20-conclusion}
We have specified a flexible, scalable benchmark that can be used both based on
pre-generated matrices or based on a procedural inclusion into an existing
finite element setting. The new model has been added to the benchmark collection
hosted at the MOR Wiki~\cite[Thermal Block]{morWiki}.

\section*{Acknowledgements}
The authors would like to thank Christian Himpe, Petar Mlinari\'c and Steffen
W.~R. Werner for helpful comments and discussions during the creation of the
model. 

Funded by the Deutsche Forschungsgemeinschaft (DFG, German Research Foundation)
under Germany's Excellence Strategy EXC 2044 –390685587, Mathematics Münster:
Dynamics–Geometry–Structure.  Funded by German Bundesministerium für Bildung und
Forschung (BMBF, Federal Ministry of Education and Research) under grant number
05M18PMA in the programme `Mathematik für Innovationen in Industrie und
Dienstleistungen'.



\end{document}